\documentclass[10pt]{amsart}

\usepackage{float}
\usepackage{amscd,amsmath,amssymb,amsfonts,amsthm, ascmac, bm}
\usepackage{enumerate, comment, mathtools} 
\usepackage{lscape}

\usepackage{caption}
\usepackage{subcaption}
\usepackage{hyperref}
\usepackage{relsize}
\usepackage{longtable}

\usepackage{booktabs}
\usepackage{a4wide}
\usepackage{verbatim}

\usepackage[mathscr]{euscript}

\usepackage{ragged2e}

\usepackage{kotex}

\usepackage{longtable,booktabs,array,tabularx,makecell}

\usepackage{xcolor}
\usepackage[linesnumbered,ruled,vlined]{algorithm2e}

\SetCommentSty{mycommfont}

\SetKwInput{KwInput}{Input}                
\SetKwInput{KwOutput}{Output}              

\usepackage{tikz-cd}
\usepackage{tikz, tikz-3dplot}
\usetikzlibrary{calc}

\numberwithin{equation}{section}
\allowdisplaybreaks[1]

\usepackage{colonequals} 

\allowdisplaybreaks

\DeclareMathOperator{\SL}{SL}
\DeclareMathOperator{\GL}{GL}

\newcommand{\C}{{\mathbb{C}}}
\newcommand{\R}{{\mathbb{R}}}

\newcommand{\Z}{{\mathbb{Z}}}

\DeclareMathOperator{\flag}{Fl}

\theoremstyle{plain}
\newtheorem{theorem}{Theorem}[section]

\newtheorem{proposition}[theorem]{Proposition}
\newtheorem{corollary}[theorem]{Corollary}

\theoremstyle{definition}

\newtheorem{definition}[theorem]{Definition}

\newtheorem{Problem}[theorem]{Problem}

\theoremstyle{remark}
\newtheorem{remark}[theorem]{Remark}

\hypersetup{colorlinks}
\hypersetup{
	unicode=false,          
	pdftoolbar=true,        
	pdfmenubar=true,        
	pdffitwindow=false,     
	pdfstartview={FitH},    
	pdftitle={Cluster-theoretic Newton--Okounkov polytopes},    
	pdfauthor={Shin Hong, Jeongsoo Kim, and Eunjeong Lee},     
	pdfnewwindow=true,      
	colorlinks=true,       
	linkcolor=blue,          
	citecolor=red,        
	filecolor=violet,      
	urlcolor=violet       
}

\begin{document}

\author{Shin Hong}
\address[Shin Hong]{School of Computer Science, 
	Chungbuk National University,
	Cheongju 28644, Republic of Korea}
\email{hongshin@chungbuk.ac.kr}

\author{Jeongsoo Kim}
\address[Jeongsoo Kim]{Department of Mathematics,
	Chungbuk National University,
	Cheongju 28644, Republic of Korea}
\email{jeongsoo.kim@chungbuk.ac.kr}

\author{Eunjeong Lee}
\address[Eunjeong Lee]{Department of Mathematics,
	Chungbuk National University,
	Cheongju 28644, Republic of Korea}
\email{eunjeong.lee@chungbuk.ac.kr}

\thanks{This work was supported by the National Research Foundation of Korea (NRF) grant funded by the Korea government (MSIT) (No.\ RS-2023-00239947).}

\title[Cluster-theoretic Newton--Okounkov polytopes]
{Cluster-theoretic Newton--Okounkov polytopes
	of \\ flag varieties of types $A_4$ and $B_2$}

\begin{abstract}
Newton--Okounkov bodies have been used for studying toric degenerations of flag varieties.	
We classify the Newton--Okounkov polytopes associated with
cluster valuations and the anticanonical line bundle on
full flag varieties of types $A_4$ and $B_2$.
In type $A_4$, the $672$ seeds of cluster type $D_6$ yield
$162$ affine unimodular equivalence classes, which coincide
with the combinatorial equivalence classes. 
We determine the subgroup of exchange-graph automorphisms
preserving every polytope class.
It is a Klein four group of index six in the automorphism
group of the exchange graph of cluster type $D_6$.
In type $B_2$, we give defining inequalities for all six seeds. 
At the anticanonical weight, the two combinatorial equivalence
classes split into three affine unimodular equivalence classes.
Their respective class-preserving automorphism groups have
orders four and two.
\end{abstract}

\date{\today}



\maketitle

\section{Introduction}

Let $G$ be a connected, simply connected semisimple algebraic
group over $\C$, and let $B\subset G$ be a Borel subgroup.
The flag variety $G/B$ is a smooth projective homogeneous variety and serves as a central object providing fruitful connections among geometry, representation theory, and combinatorics; see, e.g., \cite{Brion_Lecture,Fulton97Young}.

Although $G/B$ is not a toric variety in general, it admits various toric degenerations. 
One systematic approach is provided by the theory of
Newton--Okounkov bodies.
Under suitable finite-generation hypotheses on the associated
valuation semigroup, this construction yields toric degenerations
\cite{An13,KK12,LM09}.
On the other hand, cluster structures on double Bruhat cells $G^{u,v}$ \cite{BFZ05} (and their restriction to open unipotent cells) give a distinguished family of valuations.
Fujita--Oya \cite{FujitaOya} introduced the associated Newton--Okounkov bodies and showed that they are rational polytopes and that this family contains several known models, including Berenstein--Littelmann--Zelevinsky's string polytopes~\cite{Littelmann98} and Nakashima--Zelevinsky polytopes~\cite{NZ97}. 

Tropicalized cluster mutations relate the polytopes
associated with different seeds by piecewise-linear maps (see~\cite{FujitaOya}).
This leads to three classification questions: 
\begin{Problem}
When do two seeds give combinatorially equivalent or
affine unimodularly equivalent polytopes, which classes
are represented by string polytopes, and which symmetries
of the exchange graph preserve these classes?
\end{Problem}

In this paper, we address these questions for the toric degenerations of $G/B$ arising from the cluster structure on the double Bruhat cell $G^{e,w_0}$, where $w_0$ denotes the longest element of the Weyl group of $G$.
For a simple group $G$ of rank at least two,
Berenstein--Fomin--Zelevinsky showed that the cluster structure
on $G^{e,w_0}$ is of finite type precisely for Cartan--Killing
types $A_2$, $A_3$, $A_4$, and $B_2$, with corresponding cluster
types $A_1$, $A_3$, $D_6$, and $B_2$, respectively
\cite[Proposition~2.26]{BFZ05}.
We focus on the Cartan types $A_4$ and $B_2$, whose
corresponding cluster types are $D_6$ and $B_2$.
The former has skew-symmetric exchange matrices, whereas
the latter requires the skew-symmetrizable setting.
Table~\ref{table_finite_type} summarizes the finite type cases.
Throughout, seeds are counted up to simultaneous relabeling
of mutable variables, with the frozen variables fixed,
as specified in Section~\ref{section_NOBYs_flag5}.
A cluster structure is of \emph{finite type}  if its exchange graph
has finitely many vertices.

Our main results concern the full flag varieties of
types $A_4$ and $B_2$. We use $G=\SL_5(\C)$ in type $A_4$
and the realization $G=\mathrm{Sp}_4(\C)$ in the
isomorphic rank two setting $B_2\simeq C_2$.
For the anticanonical weight $2\rho$, set
\[
\Delta_{\bm s}
:=\Delta(G/B,\mathcal L_{2\rho},v_{\bm s},\tau_{2\rho}),
\]
where $\tau_{2\rho}$ is the distinguished section recalled
in Theorem~\ref{thm_Fujita_Oya}.
We compute these polytopes using tropicalized cluster
mutations and exact polyhedral computations in SageMath
\cite{FujitaOya,sagemath}.

Cho--Fujita--Higashitani--Lee \cite{CFHL} classified the
anticanonical cluster Newton--Okounkov polytopes in type $A_3$
and determined their class-preserving exchange-graph symmetries.
We extend this study to type $A_4$, the largest rank of finite
cluster type in the type $A$ full flag family
\cite[Proposition~2.26]{BFZ05}.
The $672$ seeds yield $162$ polytope classes, of which precisely
$23$ are represented by string polytopes at weight $2\rho$.
Together with the known cases of types $A_2$ and $A_3$, this
completes the anticanonical classification in the finite 
type $A$ setting of rank at least two.

\begin{table}
	\centering
	\begin{tabular}{l|rrrr}
		\toprule
		Cartan type of $G$ & $A_2$ & $A_3$ & $A_4$ & $B_2$ \\
		\midrule
		Cluster type & $A_1$ & $A_3$ & $D_6$ & $B_2$ \\
		Simply-laced cluster type? & yes & yes & yes & no \\
		\# of seeds & $2$ & $14$ & $672$ & $6$ \\
		\# of combinatorial classes & $1$ & $5$ & $162$ & $2$ \\
		\# of affine unimodular classes & $1$ & $5$ & $162$ & $3$ \\
		\bottomrule
	\end{tabular}
	\caption{Finite-type cases for $G^{e,w_0}$ for simple
		groups of rank at least two.
		The polytope classifications are for the anticanonical
		weight $2\rho$. See \cite{BFZ05} for the finite type
		classification, \cite{CKLP21,CFHL} for the smaller
		type $A$ cases, and
		Propositions~\ref{prop_unimodular}
		and~\ref{prop_B2_classification} for types $A_4$ and $B_2$.}
	\label{table_finite_type}
\end{table}

In type $A_4$, all $672$ anticanonical polytopes are lattice
polytopes, and translating their unique interior lattice
points to the origin makes them reflexive.
They form exactly $162$ classes under both combinatorial
and affine unimodular equivalence, and the two partitions
coincide (see Propositions~\ref{prop_lattice}
and~\ref{prop_unimodular}).
Among them, precisely $23$ classes are represented by string
polytopes for $\flag(\C^5)$ at weight $2\rho$.
The remaining $139$ classes are not even combinatorially
equivalent to a string polytope at that weight
(see Proposition~\ref{prop_string_classes}).

In type $B_2$, we give defining inequalities for all six
seeds (see Proposition~\ref{prop_B2_inequalities}).
At weight $2\rho$, all six polytopes are lattice polytopes,
but the two combinatorial classes split into three
affine unimodular classes
(see Proposition~\ref{prop_B2_classification}). 
The distinction is certified by different multisets of
edge lattice lengths.

We also determine which symmetries of the exchange graph
preserve the polytope classification at weight $2\rho$.
In type $A_4$, the exchange graph has cluster type $D_6$
and has $24$ automorphisms~\cite{CZ16finite,CZ20exchange}.
The subgroup preserving every affine unimodular equivalence
class has order four and is isomorphic to
$\mathbb Z/2\mathbb Z\times\mathbb Z/2\mathbb Z$
(see Theorem~\ref{thm_involutions}).
In type $B_2$, the corresponding subgroup has order four for
combinatorial equivalence and order two for affine unimodular
equivalence (see Corollary~\ref{cor_B2_symmetry}).
In type $A_4$, we also characterize the polytopes
$\Delta_{\bm s}$ that are affine unimodularly equivalent
to the Gelfand--Cetlin polytope for $2\rho$:
these are precisely the polytopes with $20$ facets
(see Corollary~\ref{cor_GC_type} and Problem~\ref{Prob_GC_type}).

\medskip 

The classification is established by exact computations
over $\mathbb Q$ in SageMath \cite{sagemath}.
We enumerate the finite exchange graphs, construct the
polytopes by tropicalized mutation, and compare
vertex--facet incidence graphs and lattice normal forms
of the centered polars.
For affine unimodular identifications, we recover
integral matrices of determinant $\pm1$ and integral
translation vectors, and verify their action on the
complete vertex sets.
The SageMath notebooks and verification data are available
in our GitHub repository~\cite{github}.
An accompanying interactive website\footnote{\url{https://eunjeonglee-math.github.io/flag_NOpolytopes/}}
allows readers to explore the seeds, defining inequalities,
and equivalence classes of the associated polytopes.

This paper is organized as follows.
Section~\ref{section_NOBYs_flag5} recalls the cluster
valuations and their Newton--Okounkov polytopes,
describes the computational procedures, and gives the
initial data in type $A_4$.
Section~\ref{sec_combinatorial} establishes the type
$A_4$ classification, identifies the string polytope
classes, and determines the exchange graph automorphisms.
Section~\ref{sec_B2} gives the type $B_2$ inequalities
and compares the two equivalence classifications. 
Appendix~\ref{sec_appendix} records mutation sequences,
reduced words, and $f$-vectors for the type $A_4$ classes.

\section{Description of Newton--Okounkov polytopes}\label{section_NOBYs_flag5}
\subsection{Flag varieties and their Newton--Okounkov polytopes}

In this section, we fix notation and recall the Newton--Okounkov polytopes arising from cluster valuations.
Since the general framework is standard, we keep the discussion brief and refer to
\cite{KK12,LM09,Kav15,HaradaKaveh15,FujitaOya} for background (see \cite{CFHL} for the result in type $A_3$).

Let $R$ be a domain over $\C$, and fix a total order on $\Z^m$ compatible with addition.
\begin{definition}[{\cite[Definition~3.1]{HaradaKaveh15}, \cite[Definition~2.1]{FujitaOya}}]
	A map $v\colon R\setminus\{0\}\to \Z^m$ is a \emph{valuation} if for all $\sigma,\tau\in R\setminus\{0\}$ and $c\in \C^\times$,
	\begin{enumerate}
		\item $v(\sigma\tau)=v(\sigma)+v(\tau)$,
		\item $v(c\sigma)=v(\sigma)$,
		\item $v(\sigma+\tau)\ge \min\{v(\sigma),v(\tau)\}$ whenever $\sigma+\tau\neq 0$.
	\end{enumerate}
	We say that $v$ has \emph{one-dimensional leaves} if whenever $v(\sigma)=v(\tau)$, there exists $c\in \C^\times$ such that either
	$\tau-c\sigma=0$ or $v(\tau-c\sigma)>v(\tau)$.
\end{definition}

\begin{definition}[{\cite[\S1.2]{Kav15}, \cite[Definition~2.4]{FujitaOya}}]
	Let $X$ be an irreducible normal projective variety over $\C$ of dimension $m$, and let $\mathcal L$ be a globally generated line bundle.
	Fix a valuation $v\colon \C(X)^\times\to \Z^m$ with one-dimensional leaves and a nonzero section $\tau\in H^0(X,\mathcal L)$.
	The \emph{Newton--Okounkov body} is
	\[
	\Delta(X,\mathcal L,v,\tau)
	=\overline{\bigcup_{k\ge 1}\left\{\frac{1}{k}v(\sigma/\tau^k)\ \middle|\ 0\neq \sigma\in H^0(X,\mathcal L^{\otimes k})\right\}}.
	\]
	If $\Delta(X,\mathcal L,v,\tau)$ is a polytope, we call it a \emph{Newton--Okounkov polytope}.
\end{definition}

We recall cluster-algebra terminology needed later; see \cite{FZ02_I,BFZ05} (and \cite[\S3]{FujitaOya}).
Let $J$ be a finite index set and $J_{\mathrm{uf}}\subset J$.
Set $\mathcal F=\C(z_j\mid j\in J)$.
A \emph{seed} is a pair ${\bm s}=({\bm A},\varepsilon)$ consisting of
\begin{itemize}
	\item a $J$-tuple ${\bm A}=(A_j)_{j\in J}$ of algebraically independent elements of $\mathcal F$ such that $\mathcal F = \C(A_j \mid j \in J)$, and 
	\item an integer matrix $\varepsilon=(\varepsilon_{i,j})_{i\in J_{\mathrm{uf}},\,j\in J}\in \mathrm{Mat}_{J_{\mathrm{uf}}\times J}(\Z)$
	whose $J_{\mathrm{uf}}\times J_{\mathrm{uf}}$ submatrix is skew-symmetrizable.
\end{itemize}
We allow skew-symmetrizable principal parts throughout.
In type $A_4$, the principal part is skew-symmetric and
can be represented by an ordinary quiver.
In type $B_2$, we retain the skew-symmetrizable exchange
matrix.

For $k\in J_{\mathrm{uf}}$, the \emph{mutation} $\mu_k({\bm s})=({\bm A}',\varepsilon')$ is defined by
\begin{equation}\label{eq_mutation}
	\begin{split}
		\varepsilon'_{i,j} &=
		\begin{cases}
			-\varepsilon_{i,j} & \text{if } i=k \text{ or } j=k,\\[2mm]
			\varepsilon_{i,j}+\mathrm{sgn}(\varepsilon_{i,k})\,[\varepsilon_{i,k}\varepsilon_{k,j}]_+ & \text{otherwise},
		\end{cases}\\
		A'_i &=
		\begin{cases}
			\dfrac{\prod_{j\in J}A_j^{[\varepsilon_{k,j}]_+}+\prod_{j\in J}A_j^{[-\varepsilon_{k,j}]_+}}{A_k}
			& \text{if } i=k,\\[3mm]
			A_i & \text{otherwise}.
		\end{cases}
	\end{split}
\end{equation}
Here $[a]_+=\max\{a,0\}$.

A \emph{cluster pattern} assigns a seed to each vertex of the $|J_{\mathrm{uf}}|$-regular tree $\mathbb T$, whose incident edges at each vertex are labeled bijectively by $J_{\mathrm{uf}}$, so that seeds joined by an edge labeled $k$ are related by mutation in direction $k$. 
We identify seeds under simultaneous relabeling of their mutable variables and the corresponding rows and columns of the extended exchange matrix, keeping the frozen variables fixed. The \emph{exchange graph} is the resulting quotient. 
We say that a cluster structure is \emph{of finite type} if its exchange graph has finitely many vertices, equivalently, it has finitely many seeds. 

For the type $A$ examples, take $G=\SL_n(\C)$ with
Borel subgroup $B$ and maximal torus $H\subset B$.
Then
\[
G/B \cong \flag(\C^n)
=\Bigl\{(\{0\}\subsetneq V_1\subsetneq \cdots \subsetneq V_{n-1}\subsetneq \C^n)\ \Bigm|\ \dim_\C V_i=i\Bigr\}.
\]
The Weyl group $W=N_G(H)/H\simeq \mathfrak S_n$ is generated by simple reflections $s_1,\dots,s_{n-1}$, and we write $w_0^{(n)}$ for the longest element.
Its length is $\ell(w_0^{(n)})=n(n-1)/2$.
For $n=5$, one has $\ell(w_0^{(5)})=10$, and $w_0^{(5)}$ has $768$ reduced words and $62$ commutation classes. 
Two reduced words belong to the same commutation class if they
are related by a sequence of exchanges of adjacent letters
$i,j$ satisfying $|i-j|>1$. 

We now return to the general semisimple setting.
Let $H\subset B$ be a maximal torus of $G$, let
$W=N_G(H)/H$ be its Weyl group, and let $w_0$ be its
longest element. The construction below applies to
both types $A_4$ and $B_2$.
Let $B^-\subset G$ be the opposite Borel subgroup and $U^-$ its unipotent radical.
Fix a lift $\widetilde w_0\in N_G(H)$ of $w_0$.
The \emph{open unipotent cell} is
\[
U^-_{w_0}:=U^-\cap B\widetilde w_0 B,
\]
and the standard open immersion $U^-\hookrightarrow G/B$ restricts to an open immersion $U^-_{w_0}\hookrightarrow G/B$.
The cluster structure on the double Bruhat cell $G^{e,w_0}$ from \cite{BFZ05} restricts to a cluster structure on $U^-_{w_0}$; see \cite[Appendix~B]{FujitaOya}.
We fix a cluster pattern $\mathcal S=\{{\bm s}_t\}_{t\in\mathbb T}$ for $\C[U^-_{w_0}]$.

Fujita--Oya \cite{FujitaOya} construct a valuation $v_{{\bm s}_t}$ on $\C(G/B)$ for each seed ${\bm s}_t$, and relate the corresponding Newton--Okounkov polytopes by tropicalized mutations.
We summarize the properties we will use.

\begin{theorem}[{\cite[Corollaries~5.7(2), 6.6(1), 6.30 and Theorem~6.8(1)]{FujitaOya}}]\label{thm_Fujita_Oya}
	Let $\lambda$ be a dominant integral weight and let $\mathcal L_\lambda$ be the associated line bundle on $G/B$.
	There exists a distinguished section $\tau_\lambda\in H^0(G/B,\mathcal L_\lambda)$ such that:
	\begin{enumerate}
		\item For any reduced word ${\bm i}$ of $w_0$, the Newton--Okounkov body
		$\Delta(G/B,\mathcal L_\lambda,v_{{\bm s}_{\bm i}},\tau_\lambda)$ is unimodularly equivalent to the string polytope $\Delta_{\bm i}(\lambda)$.
		\item For any vertex $t\in\mathbb T$, the Newton--Okounkov body $\Delta(G/B,\mathcal L_\lambda,v_{{\bm s}_{t}},\tau_\lambda)$ is a rational convex polytope.
		\item If $t$ and $t'$ are joined by an edge labeled $k$, then
		\[
		\Delta(G/B,\mathcal L_\lambda,v_{{\bm s}_{t'}},\tau_\lambda)
		=\mu_k^{T}\bigl(\Delta(G/B,\mathcal L_\lambda,v_{{\bm s}_{t}},\tau_\lambda)\bigr),
		\]
		where the tropicalized mutation $\mu_k^T\colon \R^{J}\to \R^{J}$ is given by
		\[
		g'_j=
		\begin{cases}
			g_j+[-\varepsilon^{(t)}_{k,j}]_+\,g_k+\varepsilon^{(t)}_{k,j}[g_k]_+ & (j\neq k),\\
			-g_k & (j=k).
		\end{cases}
		\]
		\item The string polytopes $\Delta_{\bm i}(\lambda)$ and the Nakashima--Zelevinsky polytopes $\widetilde{\Delta}_{\bm i}(\lambda)$ are related by tropicalized mutations, up to unimodular transformations.
	\end{enumerate}
\end{theorem}

The tropicalized mutation formula uses the $k$-th
mutable row of the full extended exchange matrix,
including the frozen columns.
Section~\ref{sec_B2} gives an explicit rank two example.

\subsection{Algorithms producing Newton--Okounkov polytopes}
Unless another dominant weight is explicitly specified,
we work with the anticanonical weight
\[
\lambda=2\rho=2\sum_{i=1}^{r}\varpi_i,
\qquad r=\operatorname{rank}(G),
\]
where $\varpi_1,\ldots,\varpi_r$ are the fundamental
weights. Thus $\mathcal L_{2\rho}$ is the anticanonical
line bundle on $G/B$.
For $G=\SL_n(\C)$, we have $r=n-1$; in the rank two
case of Section~\ref{sec_B2}, we have $r=2$.

In the finite type cases, the exchange graph has finitely many vertices, that is, we have  finitely many polytopes.
Starting from a chosen initial seed, we traverse the exchange graph and apply Algorithm~\ref{alg:tropical-mutation} along each mutation edge.
In practice, we implement this procedure in SageMath~\cite{sagemath}.

Algorithm~\ref{alg_tropicalized_mutation} uses the procedure of
\cite[Algorithm~1]{CFHL} to implement the tropicalized mutation
in Theorem~\ref{thm_Fujita_Oya}(3).
We split $P$ along $g_k=0$, apply the appropriate linear map
to each half, and take the convex hull of the resulting vertex images.

\begin{algorithm}[!ht]
	\DontPrintSemicolon
	\KwInput{A cluster-theoretic Newton--Okounkov polytope
		$P\subset\R^|J|$, a mutable index $k$, and the two linear maps
		$T_+$ and $T_-$ describing $\mu_k^T$ on
		$\{g_k\ge0\}$ and $\{g_k\le0\}$.}
	\KwOutput{The polytope $\mu_k^T(P)$.}
	
	$P_+ \leftarrow P\cap \{g_k\ge 0\}$; \qquad
	$P_- \leftarrow P\cap \{g_k\le 0\}$;\
	
	$V\leftarrow \emptyset$;\
	
	\For{each vertex $v$ of $P_+$}{
		$V \leftarrow V \cup \{T_+(v)\}$;
	}
	\For{each vertex $v$ of $P_-$}{
		$V \leftarrow V \cup \{T_-(v)\}$;
	}
	
	\Return{$\mathrm{conv}(V)$}\;
	\caption{Applying a tropicalized mutation to a polytope}\label{alg:tropical-mutation}
	\label{alg_tropicalized_mutation}
\end{algorithm}
By Theorem~\ref{thm_Fujita_Oya}(2)--(3),
the image $\mu_k^T(P)$ is convex in this setting.
Consequently, the convex hull returned by Algorithm~\ref{alg:tropical-mutation}
equals $\mu_k^T(P)$.
For an arbitrary input polytope, the same procedure returns
$\operatorname{conv}(\mu_k^T(P))$.

\medskip

Algorithm~\ref{alg:all-seeds} enumerates all seeds in the connected component of the initial seed in the exchange graph.
We use a standard breadth-first search (BFS): each time a new seed is discovered via a mutation, it is added to the queue and processed once.
In finite type cases this terminates after visiting all seeds.

\begin{algorithm}[!ht]
	\DontPrintSemicolon
	\KwInput{An initial seed $q_0=({\bm A}_0,\varepsilon_0)$
		with mutable index set $J_{\mathrm{uf}}$.}
	\KwOutput{The set $\mathcal S$ of all seeds reachable from $q_0$ (i.e.\ all vertices of the exchange graph component).}
	
	$\mathcal S \leftarrow \{q_0\}$;\quad
	$\mathcal Q \leftarrow$ a queue initialized with $q_0$;\
	
	\While{$\mathcal Q$ is not empty}{
			Pop $q$ from $\mathcal Q$;\
			
			\For{each $k\in J_{\mathrm{uf}}$}{
					$q' \leftarrow \mu_k(q)$;\
					
					\If{$q'\notin \mathcal S$}{
							Add $q'$ to $\mathcal S$ and push $q'$ into $\mathcal Q$;\;
						}
				}
		}
	\Return{$\mathcal S$}\;
	\caption{Enumerating all seeds by traversing the exchange graph}\label{alg:all-seeds}
\end{algorithm}
\medskip

For the remainder of this section, we illustrate the
initialization in type $A_4$, with
$G/B=\flag(\C^5)$ and $\lambda=2\rho$.
The corresponding initialization in type $B_2$ is given
in Section~\ref{sec_B2}.

Denote by $\Delta_t$ the Newton--Okounkov polytope determined by $v_t = v_{\bm s_t}$: 
\[
\Delta_t := \Delta(\flag(\C^5),\mathcal L_\lambda, v_t,\tau_\lambda). 
\]
Consider the seed~$\mathbf{s}_{t_0}$ coming from the reduced word 
\[
\bm i_0 = (1,2,1,3,2,1,4,3,2,1).
\]
The quiver associated with $\mathbf{s}_{t_0}$ is given as follows:
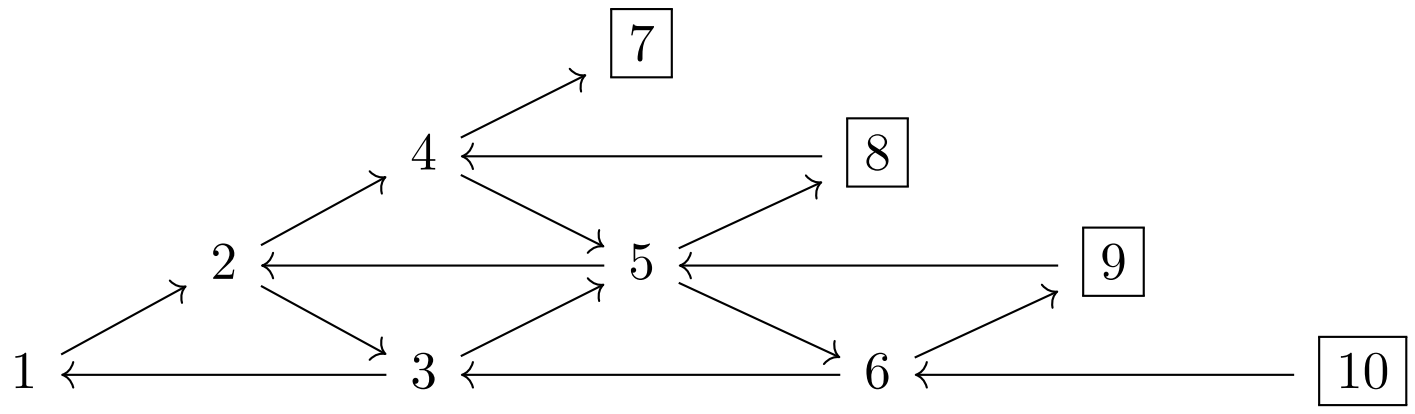
\begin{equation}\label{eq_quiver_0}
	\begin{tikzcd}[row sep = 0]
		&&& \fbox{7} \\
		&& 4 \arrow[ur] \arrow[dr] && \fbox{8} \arrow[ll] \\
		& 2 \arrow[ur] \arrow[dr] && 5 \arrow[ll] \arrow[ur] \arrow[dr] && \fbox{9} \arrow[ll] \\
		1 \arrow[ur] && 3 \arrow[ll] \arrow[ur] && 6 \arrow[ll] \arrow[ur] && \fbox{10} \arrow[ll]
	\end{tikzcd}
\end{equation}
Here, $J = \{1, \ldots, 10\}$ and $J_{\mathrm{uf}} = \{1,2,3,4,5,6\}$. 
The initial polytope has the unique interior lattice point
\[
\mathbf p=(0,0,0,0,0,0,1,1,1,1).
\]
We write $P_t:=\Delta_t-\mathbf p$ for the translated polytope.
Since the mutable coordinates of $\mathbf p$ vanish, every
tropicalized mutation satisfies
\[
\mu_k^T(\mathbf g+\mathbf p)=\mu_k^T(\mathbf g)+\mathbf p.
\]
Thus the same mutation formulas apply to the translated polytopes.
In this anticanonical setting, these polytopes are reflexive (see Proposition~\ref{prop_lattice}). In particular,
\[
P_{t_0}=\{\mathbf g\in\R^{10}\mid
\langle\mathbf g,\mathbf u_F\rangle\geq-1
\text{ for every facet }F\},
\]
where the $\mathbf u_F$ are primitive integral inward facet normals.
Here, a full-dimensional lattice polytope $P \subset \R^N$ 
with $\mathbf 0 \in \operatorname{int}(P)$ is called \emph{reflexive} if its polar $P^*$ is a lattice polytope, where $P^*$ is defined by 
\[
P^* := \{ \mathbf{y} \in \R^N \mid \langle \mathbf{y}, \mathbf{x} \rangle \geq -1 \quad \text{ for all }\mathbf{x} \in P \}.
\]

The primitive inward facet normals of the translated
initial polytope $P_{t_0}$ are the columns of the
following matrix; its defining inequalities have
constant term $1$.

\begin{equation}\label{eq:facet normal vectors}
	\resizebox{0.85\textwidth}{!}{$
		\left[
		\begin{array}{*{20}{r}}
			1 & 0 & 0 & 0 & 0 & 0 & 0 & 0 & 0 & 0 & -1 & 0 & 0 & 0 & 0 & 0 & 0 & 0 & 0 & 0 \\
			1 & 1 & 0 & 0 & 0 & 0 & 0 & 0 & 0 & 0 & 0 & -1 & 0 & 0 & 0 & 0 & 0 & 0 & 0 & 0 \\
			0 & 0 & 1 & 0 & 0 & 0 & 0 & 0 & 0 & 0 & -1 & 0 & -1 & 0 & 0 & 0 & 0 & 0 & 0 & 0 \\
			1 & 1 & 0 & 1 & 0 & 0 & 0 & 0 & 0 & 0 & 0 & 0 & 0 & -1 & 0 & 0 & 0 & 0 & 0 & 0 \\
			0&0&1&0&1&0&0&0&0&0&0&-1&0&0&-1&0&0&0&0&0\\
			0 & 0 & 0 & 0 & 0 & 1 & 0 & 0 & 0 & 0 & -1 & 0 & -1 & 0 & 0 & -1 & 0 & 0 & 0 & 0 \\
			1 & 1 & 0 & 1 & 0 & 0 & 1 & 0 & 0 & 0 & 0 & 0 & 0 & 0 & 0 & 0 & -1 & 0 & 0 & 0 \\
0&0&1&0&1&0&0&1&0&0&0&0&0&-1&0&0&0&-1&0&0\\
0&0&0&0&0&1&0&0&1&0&0&-1&0&0&-1&0&0&0&-1&0\\
0&0&0&0&0&0&0&0&0&1&-1&0&-1&0&0&-1&0&0&0&-1
		\end{array}		
		\right].
		$}
\end{equation}

These normal vectors can be obtained as follows. 
By \cite[Corollary~6.6]{FujitaOya}, the Newton--Okounkov polytope associated with $\mathbf s_{t_0}$ is related to the string polytope for ${\bm i}_0$ by an explicit unimodular linear change of coordinates, described by the matrix 
\[
M_{\bm i_0}=\begin{bmatrix}
	1&0&0&0&0&0&0&0&0&0\\
	1&1&0&0&0&0&0&0&0&0\\
	0&1&1&0&0&0&0&0&0&0\\
	1&1&0&1&0&0&0&0&0&0\\
	0&1&1&1&1&0&0&0&0&0\\
	0&0&0&1&1&1&0&0&0&0\\
	1&1&0&1&0&0&1&0&0&0\\
	0&1&1&1&1&0&1&1&0&0\\
	0&0&0&1&1&1&1&1&1&0\\
	0&0&0&0&0&0&1&1&1&1
\end{bmatrix}.
\]
We regard coordinate vectors as row vectors.

On the other hand,  using the description in~\cite[Section~5]{Littelmann98} (or  Gleizer--Postnikov's rigorous paths in~\cite{GP00}), the string polytope~$\Delta_{{\bm i}_0}(\lambda)$ is given by the set of points $(a_1,\dots,a_{10}) \in \mathbb{R}^{10}_{\geq 0}$ satisfying the following inequalities. Here, $\lambda_j=\langle\lambda,h_j\rangle$ for $1\leq j\leq4$;
in our anticanonical setting, $\lambda=2\rho$ and $\lambda_j=2$. 
\[
\begin{array}{rcccl}
	0 & \leq & a_1 & \leq & a_2 - 2a_3 + a_5 -2a_6 + a_9 - 2a_{10} + \lambda_1, \\
	a_3 & \leq & a_2 & \leq & a_3 + a_4 - 2a_5 + a_6 + a_8 -2a_9 +a_{10} + \lambda_2, \\
	0 & \leq & a_3 & \leq & a_5 - 2a_6 + a_9 -2a_{10} + \lambda_1, \\
	a_5 & \leq & a_4 & \leq & a_5 + a_7 - 2a_8 + a_9 + \lambda_3, \\
	a_6 &\leq& a_5 &\leq& a_6+a_8-2a_9+a_{10}+\lambda_2, \\
0 &\leq& a_6 &\leq& a_9-2a_{10}+\lambda_1,\\
	a_8 & \leq & a_7 & \leq & a_8 + \lambda_4, \\
	a_9 & \leq & a_8 & \leq & a_9 + \lambda_3, \\
	a_{10} & \leq & a_9 & \leq & a_{10} + \lambda_2, \\
	0 & \leq & a_{10} & \leq & \lambda_1.
\end{array}
\]
Therefore, computing
\[
\resizebox{\textwidth}{!}{$
\begin{bmatrix}
	1&0&0&0&0&0&0&0&0&0\\
	1&1&0&0&0&0&0&0&0&0\\
	0&1&1&0&0&0&0&0&0&0\\
	1&1&0&1&0&0&0&0&0&0\\
	0&1&1&1&1&0&0&0&0&0\\
	0&0&0&1&1&1&0&0&0&0\\
	1&1&0&1&0&0&1&0&0&0\\
	0&1&1&1&1&0&1&1&0&0\\
	0&0&0&1&1&1&1&1&1&0\\
	0&0&0&0&0&0&1&1&1&1
\end{bmatrix}
	\left[
	\begin{array}{*{20}{r}}
		1 & 0 & 0 & 0 & 0 & 0 & 0 & 0 & 0 & 0 & -1 & 0 & 0 & 0 & 0 & 0 & 0 & 0 & 0 & 0 \\
		0 & 1 & 0 & 0 & 0 & 0 & 0 & 0 & 0 & 0 & 1 & -1 & 0 & 0 & 0 & 0 & 0 & 0 & 0 & 0 \\
		0 & -1 & 1 & 0 & 0 & 0 & 0 & 0 & 0 & 0 & -2 & 1 & -1 & 0 & 0 & 0 & 0 & 0 & 0 & 0 \\
		0 & 0 & 0 & 1 & 0 & 0 & 0 & 0 & 0 & 0 & 0 & 1 & 0 & -1 & 0 & 0 & 0 & 0 & 0 & 0 \\
		0 & 0 & 0 & -1 & 1 & 0 & 0 & 0 & 0 & 0 & 1 & -2 & 1 & 1 & -1 & 0 & 0 & 0 & 0 & 0 \\
		0 & 0 & 0 & 0 & -1 & 1 & 0 & 0 & 0 & 0 & -2 & 1 & -2 & 0 & 1 & -1 & 0 & 0 & 0 & 0 \\
		0 & 0 & 0 & 0 & 0 & 0 & 1 & 0 & 0 & 0 & 0 & 0 & 0 & 1 & 0 & 0 & -1 & 0 & 0 & 0 \\
		0 & 0 & 0 & 0 & 0 & 0 & -1 & 1 & 0 & 0 & 0 & 1 & 0 & -2 & 1 & 0 & 1 & -1 & 0 & 0 \\
		0 & 0 & 0 & 0 & 0 & 0 & 0 & -1 & 1 & 0 & 1 & -2 & 1 & 1 & -2 & 1 & 0 & 1 & -1 & 0 \\
		0 & 0 & 0 & 0 & 0 & 0 & 0 & 0 & -1 & 1 & -2 & 1 & -2 & 0 & 1 & -2 & 0 & 0 & 1 & -1
	\end{array}
	\right],
	$}
\]
we obtain (\ref{eq:facet normal vectors}).

\section{Type $A_4$ case}
\label{sec_combinatorial}
Throughout this section, we take $G=\SL_5(\C)$ and
$\lambda=2\rho$.
We summarize the computational classification and
the symmetries on the exchange graph.

We compute the polytopes by iterating
Algorithm~\ref{alg_tropicalized_mutation} in SageMath.
The code and data are available in \cite{github}.
Table~\ref{tab:class-summary} in Appendix~\ref{sec_appendix}
lists the equivalence classes, mutation sequences, reduced words,
and $f$-vectors of representative polytopes.

We first observe that all Newton--Okounkov polytopes
considered here are lattice polytopes, that is,
all their vertices are integral vectors. Together with~\cite[Corollary~5.1]{Steinert22}, this implies that all polytopes are reflexive (after an appropriate lattice translation).

\begin{proposition}\label{prop_lattice}
For $\lambda=2\rho$,  the Newton--Okounkov polytopes 
associated with
the $672$ seeds of $\flag(\mathbb C^5)$ are lattice polytopes.
After translating their unique interior lattice points to the origin,
they are reflexive. Their associated toric varieties are
Gorenstein Fano. 
\end{proposition}
\begin{proof}
Exact computations over $\mathbb Q$ verify that all $672$
polytopes have integral vertices. The cluster valuations have
full rank, and their valuation semigroups are finitely generated
and saturated by \cite[Corollary~6.9(2)]{FujitaOya}.
Since $\lambda=2\rho$ is anticanonical,
\cite[Corollary~5.1]{Steinert22} implies reflexivity after a
lattice translation. The associated normal toric varieties are
Gorenstein Fano by \cite[Proposition~2.2.23]{Batyrev94}.	
\end{proof}

Next, we compare the polytopes up to combinatorial and unimodular equivalence using SageMath.
We use the SageMath method \verb|is_combinatorially_isomorphic()| to test combinatorial equivalence, i.e., isomorphism of face posets.
Moreover, two lattice polytopes $P, Q \subset \R^N$ are \emph{affine unimodularly equivalent} if there exist $A \in \GL_N(\Z)$ and $\mathbf{v} \in \Z^N$ such that the affine map $\mathbf{x} \mapsto A \mathbf{x} + \mathbf{v}$ sends $P$ to $Q$.
To test affine unimodular equivalence, we first translate each
polytope so that its unique interior lattice point is the origin
and then compare canonical lattice normal forms of the polars.
Since an affine unimodular equivalence maps the unique interior
lattice point to the unique interior lattice point, this reduces
the problem to linear unimodular equivalence. This yields the
following classification.

\begin{proposition}\label{prop_unimodular}
For $G=\SL_5(\C)$ and $\lambda=2\rho$, the polytopes
associated with the $672$ seeds form $162$ affine
unimodular equivalence classes.
These coincide with their combinatorial equivalence
classes \textup{(}see Table~\ref{tab:class-summary}\textup{)}.	
\end{proposition}

As presented in Table~\ref{tab:class-summary}, there are four classes of size $2$, $150$ classes of size $4$,
and eight classes of size $8$.
The number of vertices ranges from $327$ to $463$, and the number of facets ranges from $20$ to $35$.
For each class, we record representative seeds and the corresponding $f$-vectors.

Note that the $768$ reduced words form $62$ commutation classes. Commutation 
and the Dynkin involution $i\mapsto5-i$ reduce the comparison to
the $31$ representatives in \cite[Table~6]{BLMM17}. 
Comparing with the string polytopes, we obtain the following:
\begin{proposition}\label{prop_string_classes}
	For $\lambda=2\rho$, the string polytopes
	$\Delta_{\bm i}(\lambda)$ of $\flag(\mathbb C^5)$ form exactly
	$23$ classes under either combinatorial or affine unimodular
	equivalence, and these two partitions coincide.
	At $\lambda=2\rho$, precisely $23$ of the $162$ cluster polytope
	classes are represented by string polytopes. The remaining $139$
	classes are not combinatorially equivalent to any
	$\Delta_{\bm i}(\lambda)$.
\end{proposition}
\begin{proof}
	Commutation and the Dynkin involution reduce the $768$ reduced
	words to the $31$ representatives in \cite[Table~6]{BLMM17}.
	Exact lattice normal forms and vertex--facet incidence graphs
	give $23$ classes. Explicit affine unimodular maps identify these
	with $23$ distinct cluster polytope classes. The completeness of
	the representative list and the equivalence certificates are
	recorded in the supplementary material~\cite{github}.
\end{proof}

\begin{remark}
	The number $28$ reported in \cite[Proposition~4.8 and Table~6]{BLMM17}
	is reduced to $23$ by the following additional affine unimodular
	equivalences, using the labels of that table:
	\[
	S_2\cong S_{19},\quad S_3\cong S_{20},\quad
	S_4\cong S_{17},\quad S_5\cong S_{31},\quad
	S_{27}\cong S_{28}.
	\]
	For each pair, the supplementary data give an integral
	matrix $A$ with determinant $\pm1$ and an integral vector
	$\mathbf b$ such that the affine map
	$\mathbf x\mapsto A\mathbf x+\mathbf b$
	maps the complete vertex set of the first polytope
	onto that of the second.
	The transformations are verified at weight $\rho$ (and hence by scaling at $2\rho$).  
\end{remark}

The polytopes in Classes~2 and~5 have the same $f$-vector but
are not combinatorially equivalent, as their vertex-degree
distributions differ (see Table~\ref{table_degrees_Case25}).
In particular, they are not affine unimodularly equivalent.
\begin{table}
	\begin{tabular}{c|rrrrrrrrrrrr}
		\toprule 
		degree & 10 & 11 & 12 & 13 & 14 & 15 & 16 & 17 & 18 & 19 & 20 & 21 \\
		\midrule 
		Class 2 & 100 & 137 & 34 & 49 & 8 & 8 & 27 & 1 & 3 & 0 & 0 & 1 
		\\
		Class 5  &  110 & 130 & 18 & 60 & 18 & 0 & 26 & 2 & 0 & 0 & 4 & 0\\
		\bottomrule 
	\end{tabular}
	\caption{Numbers of vertices of polytopes in Classes 2 and 5  having a given degree}
	\label{table_degrees_Case25}
\end{table}

Let $\Gamma$ be the exchange graph in cluster type $D_6$. 
The results of Chang--Zhu imply
\[
\operatorname{Aut}(\Gamma)
\cong \operatorname{Aut}(\mathcal A_{D_6})
\cong \operatorname{Dih}_6\times\mathbb Z/2\mathbb Z,
\]
where $\operatorname{Dih}_6$ has order $12$
\cite{CZ16finite,CZ20exchange}.
The polytope assignment equips this graph with additional
data. We determine the subgroup that preserves each of its
affine unimodular equivalence classes.

\begin{theorem}\label{thm_involutions}
Let $\Gamma$ be the exchange graph of cluster type $D_6$
associated with $G=\SL_5(\mathbb C)$, modulo mutable-variable
relabeling and with the edge labels forgotten.
Let $c(\bm s)=[P_{\bm s}]$ be the affine unimodular equivalence
class of the anticanonical polytope associated with $\bm s$.
Then
\[
H:=\{g\in\operatorname{Aut}(\Gamma)\mid
c(g(\bm s))=c(\bm s)\text{ for every }\bm s\}
\cong\mathbb Z/2\mathbb Z\times\mathbb Z/2\mathbb Z.
\]
In particular, $[\operatorname{Aut}(\Gamma):H]=6$.
The same subgroup is obtained using combinatorial
equivalence classes.
\end{theorem}
\begin{proof} 
The exchange graph $\Gamma$ in cluster type $D_6$ has $672$ vertices and $2016$
undirected edges. Computing its full automorphism group in
SageMath gives $24$ elements, in agreement with the
type $D_6$ description above.
Testing each permutation against the class assignment of
Proposition~\ref{prop_unimodular} leaves exactly four elements
$g_0,g_1,g_2,g_3$, with $g_0$ the identity.
Their permutations satisfy
\[
g_1^2=g_2^2=1,\qquad g_1g_2=g_2g_1=g_3,
\qquad g_1\ne1,\quad g_2\ne1,\quad g_1\ne g_2.
\]
Thus $H=\langle g_1,g_2\rangle$ is a Klein four group,
and its index is $24/4=6$.
Finally, Proposition~\ref{prop_unimodular} identifies the
combinatorial and affine unimodular partitions, so the result follows.	
\end{proof}
\begin{remark}
	The action of the class-preserving subgroup $H$ on the
	$672$ seeds has $176$ orbits, whereas the polytope
	classification has $162$ classes.
	Every $H$-orbit lies in a single equivalence class, so
	the orbit partition strictly refines the polytope
	classification.
	The complete orbit list is included in the accompanying
	verification data \cite{github}.
\end{remark}

Using Proposition~\ref{prop_unimodular}, we obtain the following characterization. 
\begin{corollary}\label{cor_GC_type}
Let $P$ be a cluster-theoretic Newton--Okounkov polytope of
$\flag(\mathbb C^5)$ for $\lambda=2\rho$.
Then $P$ has $20$ facets if and only if it is affine unimodularly
equivalent to the Gelfand--Cetlin polytope for $2\rho$.
These are precisely the polytopes in Class~1 of
Table~\ref{tab:class-summary}; every other class has more than
$20$ facets.
\end{corollary}

Corollary~\ref{cor_GC_type} gives the type $A_4$ counterpart of
\cite[Corollary~3.4]{CFHL} and is consistent with the
string-polytope criterion of \cite[Theorem~A]{CKLP21}.
This suggests the following question for type $A$.
\begin{Problem}\label{Prob_GC_type}
	For $n\geq3$, let $P$ be a cluster-theoretic Newton--Okounkov
	polytope of $\flag(\C^n)$ for $\lambda=2\rho$.
	If $P$ has exactly $n(n-1)$ facets, must it be affine
	unimodularly equivalent to the Gelfand--Cetlin polytope for $2\rho$?
\end{Problem}

\section{Type $B_2$ case}\label{sec_B2}
We use the realization $G=\mathrm{Sp}_4(\mathbb C)$,
whose root system is of type $C_2\simeq B_2$.
We follow the numbering of
\cite[Example~6.33]{FujitaOya}, in which
$\alpha_1$ is the long simple root and
$\alpha_2$ is the short simple root.
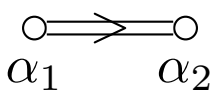
\begin{center}
	\begin{tikzpicture}
		\node[circle, draw, inner sep = 0cm, minimum size = 1.5mm, label=-90:{$\alpha_1$}] (1) at (1,0) {};
		\node[circle, draw, inner sep = 0cm, minimum size = 1.5mm, label=-90:{$\alpha_2$}] (2) at (2,0) {};

	\draw[double,  double distance=2pt] (1)--(2);
	\node at (1.5,0) {$>$};
	\end{tikzpicture}
\end{center}
The initial reduced word in this convention is
$(1,2,1,2)$. In the $B_2$ convention used in
\cite[Remark~8.5]{FujitaOya}, the same initial seed is
described by $(2,1,2,1)$ after interchanging the labels of 
the two simple roots. 
Write
$\lambda=\lambda_1\varpi_1+\lambda_2\varpi_2$,
where $\lambda_1,\lambda_2\in\mathbb Z_{\geq0}$ and
the fundamental weights use the same numbering.
The anticanonical weight is given by
$\lambda_1=\lambda_2=2$.

\begin{figure}
	\begin{tikzpicture}
		\draw[dashed] (-2.3,0)--(2.3,0) node[pos=0.95] {$\updownarrow$} node[right] {~~$\sigma$};
\draw[dashed] (0,-2.3)--(0,2.3) node[pos=0.95] {$\leftrightarrow$} node[above] {$\tau$};

		\foreach \x/\y in {0/green!50!black,1/orange,2/purple,3/purple,4/orange,5/green!50!black}{
			\node[circle, draw,fill=\y!20!white] (\x) at (-\x*60 + 150:1.5) {\x};
		}
		\foreach \x/\y in {0/1,1/2,2/3,3/4,4/5,5/0}{
			\draw (\x) -- (\y);
		}
	\end{tikzpicture}
\caption{The type $B_2$ exchange graph, colored by affine
	unimodular equivalence classes.
	The reflection $\sigma$ preserves every color;
	$\tau$ exchanges the two colors contained in the
	combinatorial class $\{0,2,3,5\}$.}\label{fig_B2}
\end{figure}
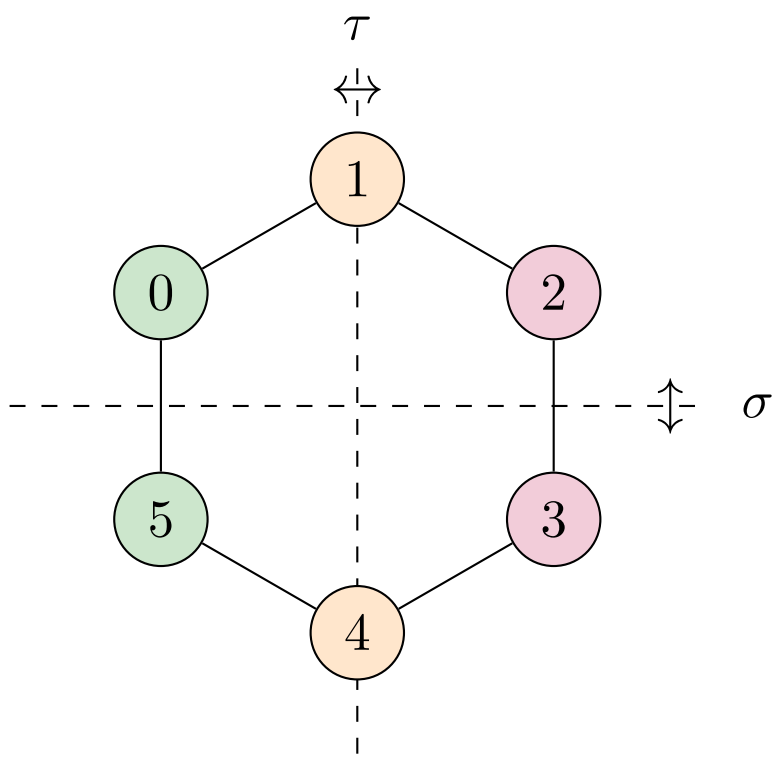

The initial extended exchange matrix is
\[
\varepsilon_0=
\begin{bmatrix}0&-2&1&0\\1&0&-1&1\end{bmatrix}.
\]
Rows correspond to the mutable indices $1,2$, and columns to
all four indices; $3,4$ are frozen. We use the tropicalized
mutation convention of Theorem~\ref{thm_Fujita_Oya}.
Starting at seed $0$, define seed $j$, for $1\leq j\leq5$,
by applying the first $j$ entries of $(1,2,1,2,1)$ in the
listed order. Mutation in direction~$2$ at seed $5$ returns
to seed $0$. The exchange graph is therefore the hexagon in Figure~\ref{fig_B2}. 

Let $\Delta_j(\lambda_1,\lambda_2)$ denote the Newton--Okounkov polytope
associated with seed $j$ and weight $\lambda_1 \varpi_1+ \lambda_2 \varpi_2$.
Here $\mathbf g=(g_1,g_2,g_3,g_4)$ denotes the unshifted
valuation coordinates. 

\begin{proposition}\label{prop_B2_inequalities}
	For every $\lambda_1, \lambda_2 \in\mathbb Z_{\geq0}$, the polytope
	$\Delta_j(\lambda_1, \lambda_2)$ is defined by
	\begin{equation}\label{eq_B2_common}
		0\leq g_3\leq \lambda_1,\qquad 0\leq g_4\leq \lambda_2,
	\end{equation}
	together with the four inequalities in row $j$ of
	Table~\ref{tab_B2_inequalities}.
	For $\lambda_1=\lambda_2=2$, all eight inequalities support facets.
\end{proposition}

\begin{table}[ht]
	\centering
	\begin{tabular}{c|ll}
		\toprule
		Seed & \multicolumn{2}{c}{Additional inequalities} \\
		\midrule
		$0$ & $g_1+g_3\leq \lambda_1$ & $g_2+g_4\leq \lambda_2$ \\
		& $g_2+g_3\geq0$ & $g_1+g_2+g_3\geq0$ \\
		\midrule
		$1$ & $2g_1+g_2+g_4\leq \lambda_2$ & $g_2+g_4\leq \lambda_2$ \\
		& $g_1+g_3\geq0$ & $g_1+g_2+g_3\geq0$ \\
		\midrule
		$2$ & $2g_1+g_2+g_4\leq \lambda_2$ & $2g_1+g_4\leq \lambda_2$ \\
		& $g_2+g_4\geq0$ & $g_1+g_3\geq0$ \\
		\midrule
		$3$ & $g_1+g_3\leq \lambda_1$ & $g_2+g_4\leq \lambda_2$ \\
		& $2g_1+g_4\geq0$ & $2g_1+g_2+g_4\geq0$ \\
		\midrule
		$4$ & $g_1+g_2+g_3\leq \lambda_1$ & $g_1+g_3\leq \lambda_1$ \\
		& $g_2+g_4\geq0$ & $2g_1+g_2+g_4\geq0$ \\
		\midrule
		$5$ & $g_1+g_2+g_3\leq \lambda_1$ & $g_2+g_3\leq \lambda_1$ \\
		& $g_2+g_4\geq0$ & $g_1+g_3\geq0$ \\
		\bottomrule
	\end{tabular}
	\caption{Defining inequalities for the six type $B_2$ seeds,
		in addition to \eqref{eq_B2_common}.}
	\label{tab_B2_inequalities}
\end{table}

The inequalities for seeds $0$ and $4$ agree with the two
cluster polytopes in \cite[Example~6.33]{FujitaOya}.
We obtain the remaining inequalities by tropicalized mutation.
We illustrate the calculation for the mutation from seed $0$ to seed $1$.
In our example
\[
\varepsilon_0^{\circ}=\begin{bmatrix}0&-2\\1&0\end{bmatrix},\qquad
D=\operatorname{diag}(1,2),\qquad
D\varepsilon_0^{\circ}=\begin{bmatrix}0&-2\\2&0\end{bmatrix}.
\]
Thus $D\varepsilon_0^{\circ}$ is skew-symmetric. 

In the notation of \cite[Equation~(4.1) and Corollary~5.7(2)]{FujitaOya},
the tropical mutation of the Newton--Okounkov body uses the entries
$\varepsilon_{k,j}$ in row $k$:
\[
h_k=-g_k,\qquad
h_j=g_j+[-\varepsilon_{k,j}]_+g_k
+\varepsilon_{k,j}[g_k]_+\quad(j\ne k).
\]

For $k=1$, the row is $(0,-2,1,0)$, so
\begin{equation}\label{eq_B2_first_tropical_mutation}
	\begin{aligned}
		h_1&=-g_1,& h_2&=g_2+2g_1-2[g_1]_+,\\
		h_3&=g_3+[g_1]_+,& h_4&=g_4.
	\end{aligned}
\end{equation}

Equation~\eqref{eq_B2_first_tropical_mutation} becomes
\[
\mathbf h=
\begin{cases}
	(-g_1,g_2,g_3+g_1,g_4),&g_1\geq0\quad(h_1\leq0),\\
	(-g_1,g_2+2g_1,g_3,g_4),&g_1\leq0\quad(h_1\geq0).
\end{cases}
\]
The corresponding inverse substitutions are
\[
\mathbf g=
\begin{cases}
	(-h_1,h_2,h_3+h_1,h_4),&h_1\leq0,\\
	(-h_1,h_2+2h_1,h_3,h_4),&h_1\geq0.
\end{cases}
\]
Substituting these expressions into every defining inequality of
$\Delta_0(\lambda_1, \lambda_2)$ gives the following list. 
\begin{center}
\begin{tabular}{r|rr}
	source inequality & $h_1\leq0 $& $h_1\geq0 $\\ \midrule 
	$g_3\geq0$ & $h_1+h_3\geq0$ & $h_3\geq0$\\
	$g_3\leq \lambda_1$ & $h_1+h_3\leq \lambda_1$ & $h_3\leq \lambda_1$ \\
	$g_4\geq0$ & $h_4\geq0$ & $h_4\geq0$ \\
	$g_4\leq \lambda_2$ & $h_4\leq \lambda_2$ & $h_4\leq \lambda_2$ \\
	$g_1+g_3\leq \lambda_1$ & $h_3\leq \lambda_1$ & $h_3-h_1\leq \lambda_1$ \\
	$g_2+g_4\leq \lambda_2$ & $h_2+h_4\leq \lambda_2$ & $2h_1+h_2+h_4\leq \lambda_2$ \\
	$g_2+g_3\geq0$ & $h_1+h_2+h_3\geq0$ & $2h_1+h_2+h_3\geq0$ \\
	$g_1+g_2+g_3\geq0$ & $h_2+h_3\geq0$ & $h_1+h_2+h_3\geq0$ 
\end{tabular}
\end{center}

Let $Q(\lambda_1, \lambda_2)$ be defined by
\begin{equation}\label{eq_B2_first_mutated_inequalities}
	\begin{gathered}
		0\leq h_3\leq \lambda_1,\qquad 0\leq h_4\leq \lambda_2,\\
		2h_1+h_2+h_4\leq \lambda_2,\qquad h_2+h_4\leq \lambda_2,\\
		h_1+h_3\geq0,\qquad h_1+h_2+h_3\geq0.
	\end{gathered}
\end{equation}
On $h_1\leq0$, the middle column above implies $h_3\geq0$
and $2h_1+h_2+h_4\leq \lambda_2$. Conversely, the inequalities of $Q(\lambda_1, \lambda_2)$
imply $h_1+h_3\leq \lambda_1$ and $h_2+h_3\geq0$ on this halfspace.
Thus, together with $h_1\leq0$, the inequalities in the
middle column define
$Q(\lambda_1,\lambda_2)\cap\{h_1\leq0\}$.
On $h_1\geq0$, the right-hand column implies
$h_1+h_3\geq0$ and $h_2+h_4\leq \lambda_2$.
Conversely, $Q(\lambda_1, \lambda_2)$ implies $h_3-h_1\leq \lambda_1$ and
$2h_1+h_2+h_3\geq0$ on this halfspace.
Consequently,
\begin{align*}
	\mu_1^T\bigl(\Delta_0(\lambda_1, \lambda_2)\cap\{g_1\geq0\}\bigr)
	&=Q(\lambda_1, \lambda_2)\cap\{h_1\leq0\},\\
	\mu_1^T\bigl(\Delta_0(\lambda_1, \lambda_2)\cap\{g_1\leq0\}\bigr)
	&=Q(\lambda_1, \lambda_2)\cap\{h_1\geq0\}.
\end{align*}
Their union equals $Q(\lambda_1, \lambda_2)$.
By \cite[Corollary~5.7(2)]{FujitaOya}, this is
$\Delta_1(\lambda_1,\lambda_2)$. Renaming $\mathbf h$ as $\mathbf g$ in
\eqref{eq_B2_first_mutated_inequalities} gives precisely row $1$
of Table~\ref{tab_B2_inequalities} and \eqref{eq_B2_common}. 
Continuing this process, we obtain the inequalities in Table~\ref{tab_B2_inequalities}.

\begin{remark}\label{rmk_B2}
With the seed indexing above, the six extended exchange matrices are
\[
\begin{aligned}
	\varepsilon_0&=
	\begin{bmatrix}0&-2&1&0\\1&0&-1&1\end{bmatrix},
	&
	\varepsilon_5&=
	\begin{bmatrix}0&2&-1&0\\-1&0&1&-1\end{bmatrix},
	\\[4pt]
	\varepsilon_1&=
	\begin{bmatrix}0&2&-1&0\\-1&0&0&1\end{bmatrix},
	&
	\varepsilon_4&=
	\begin{bmatrix}0&-2&1&0\\1&0&0&-1\end{bmatrix},
	\\[4pt]
	\varepsilon_2&=
	\begin{bmatrix}0&-2&-1&2\\1&0&0&-1\end{bmatrix},
	&
	\varepsilon_3&=
	\begin{bmatrix}0&2&1&-2\\-1&0&0&1\end{bmatrix}.
\end{aligned}
\]
In particular,
\[
\varepsilon_{5-j}=-\varepsilon_j
\qquad(0\leq j\leq5).
\]
Thus the pairing \(j\leftrightarrow5-j\) negates the full extended
exchange matrix, including its frozen columns.
\end{remark}

Set $\mathbf p=(0,0,1,1)$ and
\[
P_j:=\Delta_j(2,2)-\mathbf p\qquad(0\leq j\leq5).
\]
The point $\mathbf p$ is fixed by every tropicalized mutation,
since its mutable coordinates vanish. All the $P_j$ have
primitive facet inequalities with constant term $1$.

\begin{proposition}\label{prop_B2_classification}
	The six anticanonical polytopes are lattice polytopes, and
	their translates $P_j$ are reflexive.
	Their combinatorial equivalence classes are $	\{0,2,3,5\}, \{1,4\}$,
	whereas their affine unimodular equivalence classes are
	$\{0,5\}, \{1,4\}, \{2,3\}$.
	Their extended $f$-vectors are
	\[
	f(P_j)=
	\begin{cases}
		(1,12,26,22,8,1),&j\in\{0,2,3,5\},\\
		(1,11,24,21,8,1),&j\in\{1,4\}.
	\end{cases}
	\]
\end{proposition}

\begin{proof}
Direct computation shows that all six polytopes have integral vertices.
Their centered primitive facet inequalities have constant term $1$,
so their polars are lattice polytopes. Hence the $P_j$ are reflexive.
	
	Considering the
	inequalities gives
	\[
	P_5=A_{05}P_0,\qquad P_4=-P_1,\qquad P_3=-P_2,
	\qquad
	A_{05}=\begin{bmatrix}
		1&1&2&0\\
		0&-1&0&0\\
		0&0&-1&0\\
		0&0&0&-1
	\end{bmatrix}.
	\]
	Since $\det A_{05}=-1$, these are unimodular equivalences.
	The different $f$-vectors separate the class of $P_1$ from
	those of $P_0$ and $P_2$. An isomorphism of the vertex-facet
	incidence graphs verifies that $P_0$ and $P_2$ are
	combinatorially equivalent.
	
	To distinguish $P_0$ from $P_2$ up to affine unimodular equivalence,
	we compare their edge lattice lengths.
	For an edge with integral endpoints $u,v$, its lattice length is
	$\gcd(|u_1-v_1|,\ldots,|u_4-v_4|)$.
	The edge-length multiplicities are given as follows: 
	\begin{center}
	\begin{tabular}{c|rrrr} 
				lattice length&1&2&3&4\\\midrule 
		$P_0$&0&22&0&4
\end{tabular}
\qquad 
\begin{tabular}{c|rrrr}
		lattice length&1&2&3&6\\\midrule 
$P_2$&6&16&2&2
\end{tabular}
	\end{center}
	These multisets are invariant under affine unimodular maps
	but they are different. Consequently, $P_0$ and $P_2$ are not
	affine unimodularly equivalent, proving both classifications.
\end{proof}

Let $\Gamma$ be the exchange graph of the six type $B_2$
seeds.
We view its vertex labels as elements of
$\mathbb Z/6\mathbb Z$.
The automorphism group of the uncolored graph $\Gamma$ is
the dihedral group of order $12$. 
Define $\sigma(j)=5-j$ and $\tau(j)=2-j\pmod 6$. 
Let $\operatorname{Aut}_{\mathrm{comb}}(\Gamma)$
consist of the graph automorphisms $\phi$ for which
$P_{\phi(j)}$ is combinatorially equivalent to $P_j$
for every $j$.
Define $\operatorname{Aut}_{\mathrm{uni}}(\Gamma)$
analogously using affine unimodular equivalence.
Thus these groups preserve each class individually.
\begin{corollary}\label{cor_B2_symmetry}
	The groups preserving each combinatorial equivalence
	class and each affine unimodular equivalence class,
	respectively, are
\begin{align*}
	\operatorname{Aut}_{\mathrm{comb}}(\Gamma)
	&=\langle\sigma,\tau\rangle
	\cong\mathbb Z/2\mathbb Z\times\mathbb Z/2\mathbb Z,\\
	\operatorname{Aut}_{\mathrm{uni}}(\Gamma)
	&=\langle\sigma\rangle\cong\mathbb Z/2\mathbb Z.
\end{align*}
\end{corollary}

Place the vertices in the cyclic order $0,1,2,3,4,5$ on a
regular hexagon, with vertex $1$ at the top, vertex $4$ at
the bottom, and vertex $0$ at the upper left.
Then $\sigma$ is the reflection across the horizontal axis and
$\tau$ is the reflection across the vertical axis (Figure~\ref{fig_B2}).

\begin{remark}\label{rem_B2_literature}
	The three lattice types can be identified explicitly with
	known polytopes. At weight $2\rho$, seeds $0,5$ give the
	Gelfand--Cetlin-type string polytope in
	\cite[Example~6.33]{FujitaOya}; 
	and seeds $2,3$ represent the other affine unimodular
	string-polytope class obtained by the reduced decomposition $(2,1,2,1)$, whose string cone is described
	in \cite[Lemma~5.17]{CFL_BC}.
\end{remark}

The computations in types \(A_4\) and \(B_2\) suggest a relation between arrow reversal and unimodular equivalence of cluster Newton--Okounkov polytopes. This motivates the following question.
\begin{Problem}\label{Prob_equivalence}
For a full flag variety $G/B$ with $G$ a connected, simply connected
semisimple complex algebraic group, does there exist an involutive
automorphism $\kappa$ of its exchange graph such that, after
compatible relabeling of the mutable and frozen variable indices,
	\[
	\varepsilon_{\kappa(t)}=-\varepsilon_t,
	\qquad
	\Delta_{\kappa(t)}(2\rho)\cong_{\mathbb Z}\Delta_t(2\rho)
	\]
	for every seed \(t\), where \(\cong_{\mathbb Z}\) denotes affine unimodular equivalence?
\end{Problem}

\subsection*{Declaration of AI-assisted tools} 
ChatGPT (OpenAI) was used to assist with language editing and the development of an interactive visualization of data generated by the authors. All mathematical results and proofs were developed and verified by the authors. 

\appendix
\section{Mutation sequences and reduced words in type $A_4$}
\label{sec_appendix}

This appendix concerns the type $A_4$ family at
$\lambda=2\rho$.
A mutation sequence $(k_1,\ldots,k_\ell)$ is applied
in the listed order, starting from the initial seed
$\bm s_{t_0}$ associated with
\[
\bm i_0=(1,2,1,3,2,1,4,3,2,1).
\]
The empty sequence is denoted by $()$.
The entry \texttt{None} means that no reduced word is assigned to
the indicated seed in this table. 
We use the extended $f$-vector
$(f_{-1},f_0,\ldots,f_9,f_{10})$,
where $f_{-1}=f_{10}=1$.

The class numbers in Table~\ref{tab:class-summary}
refer to combinatorial equivalence classes in the
type $A_4$ family at weight $2\rho$.
By Proposition~\ref{prop_unimodular}, these coincide
with the affine unimodular equivalence classes
for this family.

{\footnotesize
	\setlength{\tabcolsep}{3pt}
\begin{longtable}{c|l|l}
		\caption{Combinatorial and affine unimodular equivalence
			classes of the type~$A_4$ cluster Newton--Okounkov
			polytopes for $\lambda=2\rho$, with mutation sequences,
			reduced words, and $f$-vectors.}\label{tab:class-summary} \\
		\toprule 
		Class &  Mutation sequences & Reduced words \\ 
		\midrule
		\midrule
		1 & (), & [4,3,4,2,3,4,1,2,3,4], [1,2,3,1,2,1,4,3,2,1], [4,3,2,1,4,3,2,4,3,4], [4,3,2,4,3,1,4,2,3,4], \\
		& (1, 2, 3, 1, 4, 5, 2, 6, 3, 1) & [4,3,4,2,3,1,2,4,3,4], [4,3,2,4,3,1,2,4,3,4], [4,3,2,4,1,3,4,2,3,4], [1,2,3,4,1,2,3,1,2,1], \\
		&  & [4,3,2,4,3,4,1,2,3,4], [4,3,2,4,1,3,2,4,3,4], [1,2,3,1,4,2,1,3,2,1], [1,2,1,3,2,1,4,3,2,1], \\
		&  & [4,3,4,2,1,3,4,2,3,4], [1,2,3,4,1,2,1,3,2,1], [1,2,1,3,4,2,3,1,2,1], [4,3,2,1,4,3,4,2,3,4], \\
		&  & [1,2,1,3,4,2,1,3,2,1], [1,2,3,1,4,2,3,1,2,1], [1,2,3,1,2,4,1,3,2,1], [4,3,4,2,3,1,4,2,3,4], \\
		&  & [1,2,1,3,2,4,3,1,2,1], [1,2,3,1,2,4,3,1,2,1], [1,2,1,3,2,4,1,3,2,1], [4,3,4,2,1,3,2,4,3,4] \\
\bottomrule
\end{longtable}}


\end{document}